\documentclass[preprint,12pt,3p]{elsarticle}
\usepackage[
singlelinecheck=false 
]{caption}
\usepackage{enumitem}
\usepackage{bm}
\usepackage{setspace}
\usepackage{makecell}
\usepackage{multirow}
\usepackage[dvipsnames,usenames]{xcolor}
\usepackage{amsmath}
\usepackage{amssymb}
\newtheorem{theorem}{Theorem}

\newtheorem{proposition}{Proposition}
\newtheorem{remark}{Remark}
\journal{\mbox{ \ }}
\begin{document}
\renewcommand{\arraystretch}{1.3}

\begin{frontmatter}

\title{Density and correlation in a random sequential adsorption model}

\author[label2]{Charles S. do Amaral}
\address[label2]{Departamento de Matem\'atica, Centro Federal de Educação Tecnológica de Minas Gerais, Belo Horizonte, Minas Gerais, Brazil}
\author[label1]{Diogo C. dos Santos\corref{cor1}}
\ead{diogo.santos@im.ufal.br}
\cortext[cor1]{Corresponding author.}

\address[label1]{Instituto de Matem\'atica, Universidade Federal de Alagoas, Maceió, Alagoas, Brazil}

\date{\today}

\begin{abstract}
We consider the random sequential adsorption process on the one-dimensional lattice with nearest-neighbor exclusion. In this model, each site $s \in \mathbb{Z}$ starts empty and we will try to occupy it in time $t_s$, where $(t_s)_{s\in\mathbb{Z}}$ is a sequence of independent random variables uniformly distributed on the interval $[0,1]$. The site will be occupied if both of its neighbors are vacant. We provide a method to calculate the density of occupied sites up to the time $t$, as well as the pair correlation function.
\end{abstract}

\end{frontmatter}

\section{Introduction}

	\indent The present study aims to investigate properties of a  \textit{Random Sequential Adsorption} (RSA) model. In these stochastic processes, particles are irreversibly deposited on a surface, represented by a graph, whereby the presence of particles at a particular set of sites may restrict deposition of other particles at the same or nearby sites.
  Since the 1930s, RSA models have attracted the attention of physicists and mathematicians; for an overview, see Evans \cite{evans}. Some recent works on particular RSA models can be found in \cite{Aldous,Furlan, Puljiz}.
	
We consider a type of RSA on the one-dimensional lattice $\mathbb{Z}$ named \textit{RSA with nearest-neighbor exclusion} in which particles are deposited at each site and arrive if there are no particles in their first neighbors. At most one particle can occupy a site, thus, at each site, there can only be one attempt for a particle to deposit. 

In relation to this model, analytical results regarding the density of occupied sites and the pair correlation function at any moment during deposition have already been obtained using methods such as generating functions, the independence principle, or differential equations. \cite{Dickman91, principle_2, principle_1, principle_3}. In this work, we present a method to determine these functions based solely on probabilistic arguments. Recently, applying the concept developed in this study, we derived an analytical expression for the average densities of particles in the \textit{One-dimensional AB random sequential adsorption with one deposition per site} \cite{amaral_santos}.

\subsection{\textbf{Definition and main results}}
	

Consider a sequence of \textit{i.i.d.} random variables $T=(t_s)_{s \in \mathbb{Z}}$ with uniform distribution on $[0,1]$. The corresponding product measure will be denoted by $\mathbb{P}$.  We define the continuous-time stochastic process where at time $t=0$ all sites are vacant and each site $s\in\mathbb{Z}$ will be occupied in time $t_s$ if the sites $s-1$ and $s+1$ are vacant at time $t_s$. The status of the sites \textit{occupied} and \textit{vacant} will be represented by the numbers $1$ and $0$, respectively.

Let us denote by $\omega(T,t)$ the temporal-configuration of occupied sites in $\{0,1\}^\mathbb{Z}$ and for $\omega_{s}(T,t)$ the state of the site $s$ in this configuration, both at time $t$, i.e. $\omega_{s}(t)$ is the product of the indicator functions of the sets 
\begin{equation*}
\{T\in[0,1]^{\mathbb{Z}}: T_s\leq t\},\ 
\{T\in[0,1]^{\mathbb{Z}}: \omega_{s-1}(T,T_s)=0\},\
\{T\in[0,1]^{\mathbb{Z}}: \omega_{s+1}(T,T_s)=0\}.
\end{equation*}
It is not difficult to show  that $\omega (T,t)$ is well-defined for all $t\in[0,1]$ and for $\mathbb{P}$-almost all $T\in[0,1]^{\mathbb{Z}}$ using the Harris construction. 
	
By translational invariance, we have that the probability of a site $s$ being occupied until the time $t$, 
$$\phi_{s}(t):=\mathbb{P}\left(T\in[0,1]^{ \mathbb{Z}}:\ \omega_{s}(T,t)=1\right),$$
assume the same value for all $s\in \mathbb{Z}$.
Furthermore, the translations of $\mathbb{Z}$ are ergodic with respect to $\mathbb{P}$, and then by ergodic theorems, it follows that $\phi_0(t)$ coincides with the density of occupied sites  
 up to time $t$, which is usually denoted by $\rho_t$. Pedersen and Hemmer \cite{PH93}, as well as Dickman et al \cite{Dickman91}, independently determined the exact value of $\rho_t$, demonstrating the following result.
 
  \begin{theorem}\label{theo1}
Consider the \textit{RSA with nearest-neighbor exclusion} on the one-dimensional lattice. For all $t \in [0,1]$ and all $s\in\mathbb{Z}$, we have that
    \begin{equation}
        \rho_{t}=\dfrac{1-e^{-2t}}{2}.
    \end{equation}
\end{theorem}

  The particular case $t=1$ had already been obtained by Widom \cite{Widom73}. Actually, this result when $t=1$ is equivalent to demonstrating that in the \textit{discrete Rényi packing problem} \cite{hemer}, the proportion of the interval occupied by cars approaches $1-e^{-2}$ as the length of the interval tends to infinity. This result was proven by Page in 1959 \cite{page} and had previously been observed by Flory \cite{flory}. In 2015, Gerin \cite{Gerin} employed a probabilistic construction of this model to provide a proof of the same result. Our approach shares a similarity with his, although we consider a more general case as we determine the behavior of the density at any instant during the deposition. 

The approach employed by us to determine $\rho_t$ has recently been applied to derive an analytical expression for the density in a RSA model involving the deposition of two distinct types of particles onto a lattice: particles $A$ and $B$, with the constraint that different types cannot occupy neighboring sites, and there is only one deposition attempt per site \cite{amaral_santos}. At each moment in time when a particle is deposited, it is of type $A$ with probability $\alpha$ and of type $B$ with probability $\beta=1-\alpha$. We present an analytical expression for the average densities of particles of types $A$ and $B$ at all time instances, considering all possible deposition probabilities for each particle type.

Another quantity of interest in the model studied by us is the pair correlation. Given two sites $i$ and $j$, the correlation between these two sites is defined as the covariance between the state variables $\omega_i$ and $\omega_j$. When evaluated at a given time t this quantity will be denoted by $Cov_t(i,j)$, i.e., 
\begin{align*}
    Cov_t(i,j)=\mathbb{E}[\omega_{i}(t)\omega_{j}(t)]-\mathbb{E}[\omega_{i}(t)]\mathbb{E}[\omega_{j}(t)],
\end{align*}
where $\mathbb{E}$ is the expectation operator with respect to $\mathbb{P}$. Using translation invariance one more time one can see that if $s=|i-j|$ then 
$Cov_t(i,j)=Cov_t(0,s)$. We will write
$Cov_t(0,s)=C_s(t)$.

Pedersen and Hemmer \cite{PH93}, using generating functions and the \textit{independence principle} \cite{principle_2,principle_1,principle_3}, determined the values of $C_s(t)$ demonstrating the following result.
 
 \begin{theorem}\label{theo2}
Consider the \textit{RSA with nearest-neighbor exclusion} on the one-dimensional lattice. For all $t \in [0,1]$ and all $s\in\mathbb{Z}$, the following equalities are hold
\begin{equation}
    C_s(t)=-\displaystyle\tfrac{1}{2}e^{-2t}\sum_{n=0}^{\infty}\frac{(-2t)^{2n+s+1}}{(2n+s+1)!}.
\end{equation}
\end{theorem}

\noindent The value of pair correlation for $t=1$ has already been provided by Monthus \cite{Monthus}. 

\begin{remark}
The functions obtained in Theorems \ref{theo1} and \ref{theo2} differ from those found in \cite{PH93} only because we considered times with uniform distribution instead exponential distribution. Performing an obvious change of variables one obtains that the results are the same.
\end{remark}


The goal of this work is to present proofs of Theorems \ref{theo1} and \ref{theo2} using exclusively probabilistic arguments. 

The remainder of this paper is organized as follows. Section \ref{proofT1} is dedicated to presenting the proof of the Theorem \ref{theo1}. The proof of the Theorem \ref{theo2} is presented in Section \ref{proofT2}. In both cases, we divided the organization of the sections into subsections in order to organize the presentation of the results.

\section{Proof of Theorem \ref{theo1}}\label{proofT1}
	
Define $\mathcal{A}_{t}^{s}:=\{\omega_{s}(t)=1\}$ (the site $s$ will be occupied up to time $t$). We will decompose $\mathcal{A}_{t}^{s}$ as a disjoint union of events whose probabilities we can easily calculate.
	
\subsection{\textbf{Favorable events}}
  
It is clear that for $\mathcal{A}_{t}^{s}$ to occur it is necessary to have $t_s<t$. For the other hand, if $t_s<t$ and in addition $t_s<t_{s-1}$ and $t_s<t_{s+1}$, then $\mathcal{A}_{t}^{s}$ occurs. Note that there are alternative ways for the occurrence of $\mathcal{A}_{t}^{s}$; for example, even if $t_s>t_{s+1}$ the site $s$ can still be occupied, however $t_{s+1}>t_{s+2}$ is required for this. It is easy to see that  the occurrence of the event
\begin{align*}
\mathcal{G}_{1}^{s}(t):=\left\{t>t_s,\ t_{s-1}>t_s>t_{s+1}>t_{s+2},\ t_{s+3}>t_{s+2}\right\}
\end{align*}
 implies that $\mathcal{A}_{t}^{s}$ occurs and the site $s+2$ will be occupied until the time $t_{s+2}$. In general, the occurrence of the event
\begin{align*}
		\mathcal{G}_{k}^{s}(t):=\left\{t>t_s,\ t_{s-1}>t_s>t_{s+1}>\ldots>t_{s+2k-1}>t_{s+2k},\  t_{s+2k+1}>t_{s+2k}\right\}, 
\end{align*}
for $k \geq 0$, implies that $\mathcal{A}_{t}^{s}$ occurs and all sites indexed by positive numbers that have the same parity as $s$, up to $s+2k$, will be occupied. Consider the following events defined for any positive integers $k$ and $j$
\begin{align*}
        &\mathcal{H}_{k}^{s}(t):=\{t>t_s, \ t_{s+i}>t_{s+i+1}, \mbox{\ \ \textit{for all} \ \ } i=0, 1,\ldots,k-1\}, \\ 
		&\mathcal{H}_{-j}^{s}(t):=\{t>t_s, \ t_{s-i}>t_{s-i-1}, \mbox{\ \ \textit{for all} \ \ } i=0,1,\ldots,j-1\}.
\end{align*}
We will omit the time $t$ in the notation of $\mathcal{H}_{k}^{s}(t)$. Note that, for any $j,k\geq0$ and $t \in [0,1]$, the occurrence of 
\begin{equation}
\mathcal{G}_{jk}^{s}(t)= \mathcal{H}_{-2j}^{s} \cap \left\{ t_{s-2j-1}>t_{s-2j}\right\}\cap \mathcal{H}_{2k}^{s}\cap \{t_{s+2k+1}>t_{s+2k}\}
		\label{favoravel}
\end{equation}
implies that all sites between $s-2j$ and $s+2k$ that have the same parity as $s$ will be occupied up to  time $t_s$ and in particular $\mathcal{A}_{t}^{s}$ occurs. Furthermore,
\begin{equation}\label{s_occ}
\mathcal{A}_{t}^{s}=\bigcup_{j, k \geq 0} \mathcal{G}_{jk}^{s}(t) 
\end{equation} 
\noindent and the events $\mathcal{G}^s_{jk}(t)$ are mutually exclusive, whence

\begin{align}\label{phi_t}
\phi_s(t)=\sum_{j\geq0}\sum_{k\geq0} \mathbb{P}(\mathcal{G}_{jk}^{s}(t)).
\end{align}

\subsection{\textbf{Evaluation  of the probabilities}}\label{probabilidade_calculo}
 
Let us here compute the probability of the events $\mathcal{G}_{jk}^{0}(t)$, denoted by $P_{jk}(t)$. When $s=0$, we will omit the symbol $s$ in the notation of $\mathcal{H}_{k}^{s}(t)$ and $\mathcal{H}_{-j}^{s} (t)$. Note that for any two events $\mathcal{X}$ and $\mathcal{Y}$, it is hold that 
\begin{align}\label{fundamental}
\mathbb{P}(\mathcal{X}\cap\mathcal{Y})=\mathbb{P}(\mathcal{X})-\mathbb{P}(\mathcal{X}\cap\mathcal{Y}^c).
\end{align}
It follows directly from \eqref{fundamental} that
\begin{align}\label{pjk} 
P_{jk}(t)&=\mathbb{P}\left(\mathcal{H}_{-2j},\ \mathcal{H}_{2k},\ t_{2k+1}>t_{2k}\right)-\mathbb{P}\left(\mathcal{H}_{-2j-1},\ \mathcal{H}_{2k},\ t_{2k+1}>t_{2k}\right).
	\end{align}
Applying \eqref{fundamental} twice more, with $\mathcal{Y}=\{t_{2k+1}>t_{2k}\}$, we can rewrite the Equation \eqref{pjk} as
\begin{equation}\label{TiraBota}
P_{jk}(t)=\mathbb{P}\left(\mathcal{H}_{-2j},\ \mathcal{H}_{2k}\right)-\mathbb{P}\left(\mathcal{H}_{-2j},\ \mathcal{H}_{2k+1}\right)
-\mathbb{P}\left(\mathcal{H}_{-2j-1},\ \mathcal{H}_{2k}\right)
+\mathbb{P}\left(\mathcal{H}_{-2j-1},\ \mathcal{H}_{2k+1}\right). 
\end{equation}
The probabilities $\mathbb{P}(\mathcal{H}_{-m},\mathcal{H}_n)$ can be evaluated as follows
\begin{align}	
\mathbb{P}\left(\mathcal{H}_{-m},\ \mathcal{H}_{n}\right)&=\int_0^t \mathbb{P}\left(u>t_{-1}>\ldots>t_{-m},\ u>t_1>\ldots>t_{n}\right) du\nonumber \\ 
&=\int_0^t\frac{u^{n}}{n!}\frac{u^{m}}{m!} \ du=\frac{t^{m+n+1}}{m!n!(m+n+1)}.\label{probabilidadeMN}
\end{align}

\noindent Therefore, from (\ref{TiraBota}) and (\ref{probabilidadeMN}) we have that
\begin{align}
    P_{jk}(t)&=\frac{t^{2k+2j+1}}{(2k)!(2j)!(2k+2j+1)}-\frac{t^{2k+2j+2}}{(2k+1)!(2j)!(2k+2j+2)}\nonumber\\ 
&-\frac{t^{2k+2j+2}}{(2j+1)!(2k)!(2k+2j+2)}+\frac{t^{2k+2j+3}}{(2k+1)!(2j+1)!(2k+2j+3)}.\label{final_pjk}
\end{align}

We will denote the four parcels on the \textit{r.h.s.} of the above equation by $a_{jk}$, $b_{jk}$, $c_{jk}$ and $d_{jk}$, in that order, and the variable $t$ will be omitted from this notation. 

\subsection{\textbf{Calculation of series}}

Fixed $k\in\mathbb{N}$ consider
\begin{align*}
		\mathcal{I}_k(t)= \sum_{j\geq0} (a_{jk}- b_{jk}- c_{jk}+ d_{jk}).
	\end{align*}
In other words, $\mathcal{I}_k=\sum_{j\geq0}P_{jk}$. Expanding the hyperbolic sine and cosine functions as a power series, it can be easily verified that
\begin{align*}
\sum_{j\geq0}a_{jk}=\int_0^t\tfrac{u^{2k}}{(2k)!}\cosh u\ du,\ \ \ \ 
		\sum_{j\geq0} b_{jk}=\int_0^t \tfrac{u^{2k+1}}{(2k+1)!}\cosh{u} \ du, \\
\sum_{j\geq0} c_{jk} =\int_0^t\tfrac{u^{2k}}{(2k)!}\sinh{u} \ du,\ \ \mbox{and} \ \ \ 
\sum_{j\geq0}d_{jk}=\int_0^t\tfrac{u^{2k+1}}{(2k+1)!}\sinh{u}\ du.
	\end{align*} 
When we combine all together, we get that
\begin{align}\label{lindona}
\mathcal{I}_k(t)=\displaystyle\int_0^t\left[\tfrac{u^{2k}}{(2k)!} - \tfrac{u^{2k+1}}{(2k+1)!}\right]e^{-u} du=\tfrac{t^{2k+1}}{(2k+1)!}\cdot e^{-t}.
\end{align}
Therefore, Equations \eqref{s_occ} and \eqref{lindona} implies that
\begin{align}
\phi (t)=\sum_{k\geq0}\mathcal{I}_k(t)=e^{-t}\cdot\sinh{t}=\frac{1-e^{-2t}}{2}.
		\label{final_rho}
	\end{align}

\section{Proof of Theorem \ref{theo2}} \label{proofT2}

In this section, we will calculate the $C_s$ pair correlation function.
First observe that (omitting the variable $t$ time)
\begin{align*}
\phi_s&=\mathbb{P}(\omega_{-1}=\omega_{s}=1,\ \omega_{0}=0) + \mathbb{P}(\omega_{-1}=0,\ \omega_{0}=\omega_{s}=1)+
		 \mathbb{P}(\omega_{-1}=\omega_{0}=0,\ \omega_{s}=1).
\end{align*}
Then, defining $\mathcal{D}_{s}(t)=\{\omega_{-1}(t)=\omega_{0}(t)=0,\ \omega_{s}(t)=1\}$ and $\gamma_{s}(t)=\mathbb{P}(\mathcal{D}_{s}(t))$, and denoting by $p_s(t)$ the probability that both origin and $s$ are occupied until the time $t$,   we have that 
\begin{align}\label{soma_prob}
\phi_s=p_{s+1} + p_{s} + \gamma_{s}. 
\end{align}

The value of $\gamma_s$ is given by the following proposition in the case that $s$ is an even number. It will be proved in Subsection \ref{proof_theorem2}. 

\begin{proposition}\label{proposition1}
For any $t\in[0,1]$ and every $s$ even number, it is held that 
$$\gamma_s(t)=-\tfrac{1}{2} e^{-2t} \cdot \sum_{i=1}^{s} \dfrac{(-2t)^{i}}{i!}.$$
\end{proposition}

\begin{remark}
In the proof of Proposition \ref{proposition1}, the choice of $s$ being even will be important to decompose the event $\mathcal{D}_{s}(t)$ as a disjoint union of events whose probabilities can be calculated without much effort. When $s$ is odd, the same decomposition leads to event intersections whose calculations are complicated. As will be seen below, the calculation of $\gamma_s(t)$ for $s$ even is enough to determine $C_s(t)$ for all $s \in \mathbb{Z}$.
\end{remark}

We continue with the calculation of $C_s$. Suppose $r$ is an even number. Substituting the formula found in Proposition \ref{proposition1} in Equation \eqref{soma_prob} and manipulating $\phi_r-\gamma_r$ we obtain that
\begin{align}\label{SomaProbab}
    p_{r+1}(t) + p_{r}(t) & =2 \phi_{r}(t)^2- \tfrac{1}{2} e^{-2t} \sum_{i=r+1}^{\infty} \tfrac{(-2t)^{i}}{i!}.
\end{align}
Note also, for every $r$ natural number, it is holds that 
\begin{align}\label{cov}
	C_r=\mathbb{P}(\omega_{0}=\omega_{r}=1)-\phi_r^2.
\end{align}
Therefore, by Equations  \eqref{SomaProbab} and \eqref{cov}, we have that 
\begin{align}\label{cov_soma}
    C_{r+1}(t)+C_{r}(t)=- \tfrac{1}{2} e^{-2t} \sum_{i=r+1}^{\infty} \tfrac{(-2t)^{i}}{i!}.
\end{align}
For $n\in\mathbb{N}$, consider $\Gamma_{t,n}(s)=C_{s+2n+1}(t)+C_{s+2n}(t)$. Writing $C_s$ as a telescopic series and using Equation \eqref{cov_soma} twice, with $r=s+2n$ and $r=s+2n+1$, we obtain that
\begin{align}\label{sPar}
     C_{s}(t) &= \sum_{n=0}^{\infty} [\Gamma_{t,n}(s)-\Gamma_{t,n}(s+1)]
    =-\displaystyle\tfrac{1}{2}e^{-2t}\sum_{n=0}^{\infty}\tfrac{(-2t)^{2n+s+1}}{(2n+s+1)!}.
\end{align}

The formula for $C_s$ when s is an odd number can be easily obtained from Equations \eqref{cov_soma} and \eqref{sPar}.


\subsection{\textbf{Proof of Proposition \ref{proposition1}}}\label{proof_theorem2}

This subsection aims to show that, for any $t\in[0,1]$ and every $s$ even number, it is held that 
$$\mathbb{P}\left(\omega_{-1}(t)=\omega_{0}(t)=0,\ \omega_{s}(t)=1\right)=-\tfrac{1}{2} e^{-2t} \sum_{i=1}^{s} \dfrac{(-2t)^{i}}{i!}.$$
The above probability was denoted by $\gamma_s(t)$. For that purpose consider the following events which indicate whether the times $t_0$ and $t_{-1}$ are greater or less than $t$,
\begin{align*}
&\mathcal{M}_1(t):=\{t_{-1}>t,\ t_0 > t\},\ \mathcal{M}_2(t):=\{t_{-1}<t,\ t_0<t\}, \\
& \mathcal{M}_3(t):=\{t_{-1}<t< t_0\}, \
\mathcal{M}_4(t):=\{t_{-1}>t>t_0\}. \
\end{align*}
In turn, consider $\gamma_{s,i}(t):=\mathbb{P}\left(\mathcal{D}_{s}(t)\cap\mathcal{M}_i(t)\right)$ and note that 
\begin{equation}\label{gamma}
\gamma_s=\gamma_{s,1}+\gamma_{s,2}+\gamma_{s,3}+\gamma_{s,4}.    
\end{equation}

We will calculate $\gamma_{s,i}(t)$ for  $i=1$ and $i=2$, the others can be calculated similarly. For simplicity, we will omit the symbol $t$ in the notation of  $\mathcal{M}_i(t)$.  The strategy to calculate $\gamma_{s,i}$ will be to decompose $\mathcal{D}_{s}(t)\cap\mathcal{M}_i$ as unions of intersections of independent events. This is done by analyzing the random variables $(t_s)_{s \in \mathbb{Z}}$ and discarding the intersections of mutually exclusive events.

We start with the calculation of $\gamma_{s,1}$. Note that $\gamma_{s,1}(t)=\mathbb{P}(\mathcal{M}_1,\  \omega_{s}(t)=1)$. Observe also that, for all $l\geq s/2$ and $m\geq0$ it is holds that $\mathcal{G}_{lm}^s\cap\{t_0>t\}=\emptyset$, hence we have that 
\begin{align*}
\mathcal{M}_1\cap\left(\bigcup_{l,m\geq0}\mathcal{G}_{lm}^s\right)=\mathcal{M}_1\cap\left(\bigcup_{m\geq0}\bigcup_{l=0}^{(s-2)/2}\mathcal{G}_{lm}^s\right).
\end{align*}
For the other hand, if $l\leq (s-2)/2$, the events  $\mathcal{M}_1$ and $\mathcal{G}_{lm}^s$ are independent for every $m$. Therefore, using Equation \eqref{lindona} one obtain that
\begin{align*}
\gamma_{s,1}(t)=(1-t)^2 \sum_{l=0}^{(s-2)/2}\sum_{m \geq 0} \mathbb{P}(\mathcal{G}_{lm}^{s})= (1-t)^2 e^{-t} \sum_{l=0}^{(s-2)/2} \frac{t^{2l+1}}{(2l+1)!}.
\end{align*}
We will then compute $\gamma_{s,2}$. First remember the definition \eqref{favoravel} and note that 
\begin{align*}
    \gamma_{s,2}&=\mathbb{P}(\mathcal{M}_2,\ \omega_{-2}(t)=\omega_{1}(t)=\omega_{s}(t)=1)
    =\mathbb{P}\left(\mathcal{M}_2,\ \bigcup_{u, v \geq 0} \mathcal{G}_{uv}^{-2}, \ \bigcup_{j,k \geq 0} \mathcal{G}_{jk}^{1}, \ \bigcup_{l,m \geq 0} \mathcal{G}_{lm}^{s} \right).
\end{align*}
Next observe that either $k\geq(s-2)/2$ or $l\geq (s-2k-2)/2$,  then $\mathcal{G}_{jk}^1$ and $\mathcal{G}_{lm}^s$ are disjoint events. Otherwise, they are independent events. Thus,
\begin{align*}
\gamma_{s,2}&=\mathbb{P}\left(\mathcal{M}_2,\ \bigcup_{u,v, m, j \geq 0} \bigcup_{k=0}^{(s-4)/2\ \ } \bigcup_{l=0}^{(s-2k-4)/2} \mathcal{G}^{-2}_{uv} \cap \mathcal{G}^1_{jk} \cap \mathcal{G}_{lm}^{s}  \right)\\
&=\sum_{u,v,m,j\geq0} \sum_{k=0}^{(s-4)/2\ \ } \sum_{l=0}^{(s-2k-4)/2} \mathbb{P}\left(\mathcal{M}_2\cap \mathcal{G}_{uv}^{-2}\cap \mathcal{G}^1_{jk}\right)\mathbb{P}\left(\mathcal{G}_{lm}^s\right).
\end{align*}
Denoting  $\mathcal{B}_k(t):=\{t>t_0>t_1>...>t_{2k+1}, \ t_{2k+2}>t_{2k+1}\}$ and
$$\mathcal{B}_{-u}(t):=\{t>t_{-1}>t_{-2}>...>t_{-2u-2}, \ t_{-2u-3}>t_{-2u-2}\}$$ 
it is not difficult to verify that  
$\mathcal{M}_2\cap\mathcal{G}_{u0}^{-2}\cap\mathcal{G}_{0k}^1= \mathcal{B}_{-u}(t)\cap\mathcal{B}_k(t)$. For the other hand, if either $v$ or $j$ are positive we have that $\mathcal{G}_{uv}^{-2}\cap\mathcal{G}_{jk}^1=\emptyset$. Then, 
\begin{align*}
    \gamma_{s,2}&=\sum_{u,m\geq0} \sum_{k=0}^{(s-4)/2\ \ } \sum_{l=0}^{(s-2k-4)/2} \mathbb{P}\left( \mathcal{B}_{-u}\right)\mathbb{P}\left(\mathcal{B}_k\right)\mathbb{P}\left(\mathcal{G}_{lm}^s\right).
\end{align*}
With \eqref{fundamental} and \eqref{probabilidadeMN} in mind, one can easily verified that
   $\mathbb{P}(\mathcal{B}_k(t))=\mathbb{P}(\mathcal{B}_{-u}(t))=\tfrac{t^{2k+2}}{(2k+2)!}-\tfrac{t^{2k+3}}{(2k+3)!}$, for $k=u$. It is easy to see that 
		$\displaystyle\sum_{j \geq 0} \mathbb{P}(\mathcal{B}_k(t))=e^{-t}-1+t$. Thus, by using (\ref{lindona}), 
  
\begin{align}\label{gamma2}
\gamma_{s,2}(t)=(e^{-t}-1+t)e^{-t} \sum_{k=0}^{(s-4)/2} \sum_{l=0}^{\ \ (s-2k-4)/2} \left[ \frac{t^{2l+2k+3}}{(2l+1)!(2k+2)!} - \frac{t^{2l+2k+4}}{(2l+1)!(2k+3)!} \right].
\end{align}
 The double sum of the first part of \eqref{gamma2}, taking $s=2r$ for some natural number $r$, is given by
\begin{align*}
S_1(t):=\sum_{k=0}^{r-2} \sum_{l=0}^{r-k-2} \dfrac{t^{2l+2k+3}}{(2l+1)!(2k+2)!} &= \sum_{i=0}^{r-2}\sum_{j=0}^{i} \dfrac{t^{2i+3}}{(2j+1)!(2i-2j+2)!}\\ 
&=\sum_{i=0}^{r-2} \dfrac{t^{2i+3}}{(2i+3)!}\sum_{j=0}^{i} \binom{2i+3}{2j+1}\\
&= \sum_{i=0}^{r-2} \dfrac{t^{2i+3}}{(2i+3)!}\left(2^{2i+2} -1 \right)\\
&=\tfrac{1}{2}\sum_{i=0}^{r-2} \dfrac{(2t)^{2i+3}}{(2i+3)!}-\sum_{i=0}^{r-2} \dfrac{t^{2i+3}}{(2i+3)!}.
\end{align*}
Analogously, it can be verified that the double sum of the second part of Equation \eqref{gamma2} is given by
\begin{align*}
\begin{split}		
S_2(t):=&\sum_{k=0}^{r-2} \sum_{l=0}^{r-k-2} \frac{t^{2l+2k+4}}{(2l+1)!(2k+3)!} = \tfrac{1}{2}\sum_{i=0}^{r-2} \dfrac{(2t)^{2i+4}}{(2i+4)!}-\sum_{i=0}^{r-2} \dfrac{t^{2i+4}}{(2i+3)!}.
\end{split}
\end{align*}
Therefore, we can group and simplify some terms to get that
\begin{align}\label{somasoma}
S_1(t)-S_2(t)=-\tfrac{1}{2}  \sum_{i=3}^{2r} \frac{(-2t)^{i}}{i!} \ - \ (1-t) \sum_{i=0}^{r-2} \frac{t^{2i+3}}{(2i+3)!}.
\end{align}
Substituting \eqref{somasoma} in (\ref{gamma2}),
remembering that $2r=s$,
we have that
\begin{align*}
\gamma_{s,2}(t)=-\left(e^{-t}-1+t\right)e^{-t} \left[\tfrac{1}{2} \sum_{i=3}^{s} \frac{(-2t)^{i}}{i!}  + (1-t) \sum_{i=0}^{(s-4)/2} \frac{t^{2i+3}}{(2i+3)!} \right].
\end{align*}
Through reasoning analogous to those used in the calculations of $\gamma_{s,1}(t)$ and $\gamma_{s,2}(t)$, one can show that

\begin{align*}
		\gamma_{s,3}(t)&= (1-t) \left(e^{-t}-1+t\right)e^{-t} \sum_{i=0}^{(s-2)/2} \frac{t^{2i+1}}{(2i+1)!},\\
\gamma_{s,4}(t)&=-(1-t)e^{-t} \left[\tfrac{1}{2} \sum_{i=3}^{s} \frac{(-2t)^{i}}{i!} + (1-t) \sum_{i=0}^{(s-2)/4} \frac{t^{2i+3}}{(2i+3)!} \right].
\end{align*}
Next, it is easy to verify that
\begin{align}
	\gamma_{s,1}(t)+\gamma_{s,3}(t)&= e^{-2t} (1-t)  \sum_{i=0}^{(s-2)/2} \frac{t^{2i+1}}{(2i+1)!},\label{soma1}\\
	\gamma_{s,2}(t)+\gamma_{s,4}(t) &= -e^{-2t} \left[\tfrac{1}{2}  \sum_{i=3}^{s} \dfrac{(-2t)^{i}}{i!}  + (1-t) \sum_{i=0}^{(s-4)/2} \dfrac{t^{2i+3}}{(2i+3)!} \right]\label{soma2}.
\end{align}
Finally, by combining Equations \eqref{gamma}, \eqref{soma1} and \eqref{soma2} we reach the desired result,
\begin{align*}
	\gamma_{s} &=-e^{-2t} \left[\tfrac{1}{2}  \sum_{i=3}^{s} \frac{(-2t)^{i}}{i!} - (1-t)t \right] = -\tfrac{1}{2} e^{-2t}   \sum_{i=1}^{s} \frac{(-2t)^{i}}{i!}.
\end{align*}

\section*{Acknowledgements}
The authors are grateful to Mate Puljiz of the University of Zagreb for letting us know about the paper \cite{Gerin}. D.C.S. was partially supported by CNPq (Grant 409198/2021-8) and FAPEMIG (Processo APQ-00774-21).

\end{document}